\documentclass[12pt]{article}

\usepackage{graphicx} 
\usepackage{url}      
\newcommand{\doi}[1]{\url{http://dx.doi.org/#1}}
\usepackage{amsmath}  
\usepackage{amssymb}
\usepackage{mathrsfs}

\usepackage{xspace,algorithm,algorithmic,tabularx}

\usepackage{setspace}
\usepackage{color}
\usepackage{pdflscape}
\usepackage{graphics}
\usepackage{rotating}
\usepackage{verbatim}
\usepackage{amsmath}
\usepackage{amssymb}
\usepackage{amsfonts}
\usepackage{mathrsfs}
\usepackage{array}
\usepackage{epic,eepic}
\usepackage{fancyhdr}
\usepackage{url}
\usepackage{subfigure}

\newtheorem{theorem}{Theorem}[section]

\newtheorem{remark}{Remark}[section]

\newtheorem{corollary}{Corollary}[section]

\newenvironment{proof}{\begin{trivlist}
    \item[\hskip\labelsep{\it Proof.}]}{$\hfill\Box$ \end{trivlist}}

\def\bsx{\boldsymbol{x}}
\def\bsy{\boldsymbol{y}}

\def\bsv{\boldsymbol{v}}

\def\bs0{\boldsymbol{0}}

\newcommand\Var{\textnormal{Var}}

\allowdisplaybreaks

\title{Scrambled Polynomial Lattice Rules for Infinite-Dimensional Integration}

\author{Jan Baldeaux}

\begin{document}
\maketitle

\begin{abstract}

In the random case setting, scrambled polynomial lattice rules as discussed in \cite{BD10} enjoy more favourable strong tractablility properties than scrambled digital nets. This short note discusses the application of scrambled polynomial lattice rules to infinite-dimensional integration. In \cite{HMGNR10}, infinite-dimensional integration in the random case setting was examined in detail, and results based on scrambled digital nets were presented. Exploiting these improved strong tractability properties of scrambled polynomial lattice rules and making use of the analysis presented in \cite{HMGNR10}, we improve on the results that were achieved using scrambled digital nets.

\end{abstract}

\section{Introduction} \label{introduction}

In a recent series of papers, \cite{CDMGR09}, \cite{KSWW09}, \cite{HMGNR10}, \cite{NHMGR10}, the problem of infinite-dimensional quadrature has been studied. Such problems have many applications, e.g. in mathematical finance, see \cite{NH08}, \cite{IK10}, where series expansions are used to represent particular random variables.

This short note focuses on the random case setting, which was addressed in \cite{HMGNR10}. In the latter paper, a complete and general analysis was presented, clearly showing the reader how to employ a given quadrature rule. Furthermore, results based on the scrambled Niederreiter sequence were presented. Recently, \cite{BD10}, it was shown that for multivariate integration in the random case setting, scrambled polynomial lattice rules possess more favorable strong tractability properties than any scrambled digital sequences.

This raises the natural question, whether the results based on the scrambled Niederreiter sequence can be improved using scrambled polynomial lattice rules. In this short note, we give an affirmative answer to this question. In particular, the contribution of this note is the following: Using the analysis presented in \cite{HMGNR10}, we employ the multivariate in tegration result from \cite{BD10} to improve on the results presented in \cite{HMGNR10}. Furthermore, due to the growth of the $t$-value of order $t \asymp s$, it is clear that the results as presented in this paper cannot be achieved using digital sequences. Finally, we remark that results on multivariate integration using lattice rules in the random setting have not appeared in the literature.

We place ourselves in the same setting as discussed in \cite{HMGNR10}, use the same algorithms, but employ a different quadrature rule, in particular with better strong tractability properties. In the interest of giving due credit to the authors of \cite{HMGNR10} and also of saving space, we have decided to proceed as follows: We very briefly recall the function space introduced in \cite{HMGNR10} and the sampling regimes, but introduce cost and worst-case errors under the premise that algorithms are based on scrambled polynomial lattice rules. The interested reader is referred to \cite{HMGNR10} for a complete and general treatment of the problem studied in this note.

\section{The Setting} \label{sectheseeting}

In this section, we briefly recall the function space and the sampling regimes, cost and errors as introduced in \cite{HMGNR10}. Regarding notation, as in \cite{HMGNR10}, $v$ is used to denote finite subsets of $\mathbb{N}$, the set $\left\{1,\dots,s \right\}$ is denoted by $1:s$, lastly, we write $x_k \preceq y_k$ for sequences of positive real numbers $x_k$ and $y_k$, if $x_k \leq c y_k$ is valid for $k \in \mathbb{N}$ and some constant $c>0$.

\subsection{The Function Space}

We briefly remind the reader how to construct functions of infinitely many variables, as presented in \cite{HMGNR10}. Essentially, we start with a one-dimensional reproducing kernel Hilbert space, construct spaces of finitely many variables as tensor product spaces and take limits to allow for infinitely many variables. Coordinate weights $\gamma_v=\prod_{j \in v} \gamma_j$, $v \subset \mathbb{N}$, which indicate teh importance of the variables $x_j$, $j \in v$, ensure convergence of the relevant quantities.

In particular, we consider the reproducing kernel
\begin{displaymath}
k(x,y) = \frac{1}{3} + (x^2 + y^2)/2 - \max(x,y) \, ,
\end{displaymath}
$x,y \in [0,1]$, and consider the Hilbert space $H(1+\gamma k)$, for a weight $\gamma >0$, whose norm satisfies
\begin{displaymath}
\| f \|^2 = \left( \int^1_0 f(y) dy \right)^2 + \gamma^{-1} \int^1_0 (f')^2 (y) dy \, .
\end{displaymath}
To allow for functions of finitely many variables, we consider the reproducing kernel
\begin{displaymath}
K_v (\bsx, \bsy) = \prod_{j \in v} (1 +\gamma_j k(x_j,y_j))
\end{displaymath}
and of course the associated Hilbert space $H(K_v)$ is of tensor product form
\begin{equation} \label{eqfindimHspace}
H(K_v) = \bigotimes_{j \in v} H(1+\gamma_j k) \, .
\end{equation}
To define functions of infinitely many variables,w e define the measurable kernel $K$ on $[0,1]^{\mathbb{N}} \times [0,1]^{\mathbb{N}}$
\begin{displaymath}
K(\bsx, \bsy) = \sum_{v} \gamma_v K_v(\bsx, \bsy) = \sum_{v} \gamma_v \prod_{j \in v} k(x_j,y_j) \, ,
\end{displaymath}
for $\bsx, \bsy \in [0,1]^{\mathbb{N}}$ and denote the associated space by $H(K)$, which, see \cite[Lemma 6]{HMGNR10}, consists of all functions
\begin{displaymath}
f = \sum_{v} f_v \, , \quad f_v \in H_v \, ,
\end{displaymath}
for which
\begin{displaymath}
\sum_{v} \gamma^{-1}_v \| f_v \|^2_{k_v} \le \infty \, ,
\end{displaymath}
and, in case of convergence,
\begin{displaymath}
\| f \|^2_K = \sum_v \gamma^{-1}_v \| f \|^2_{k_v} \, .
\end{displaymath}

\subsection{Sampling regimes, cost, and worst-case error} \label{subsecsamplingregcostwce}

In this subsection, we introduce randomized algorithms for the integration of functions $f: [0,1]^{\mathbb{N}} \rightarrow \mathbb{R}$; the reader is referred to \cite{TWW88} and \cite{CDMGR09} for a detailed discussion.

Follwoing \cite{CDMGR09}, \cite{HMGNR10}, two sampling regimes, which specify the domains from which the integration nodes can be specified, are introduced.

{\em Fixed subspace sampling} restricts this domain to a finite-dimensional affine subspace
\begin{displaymath}
\mathcal{X}_{v,a} = \left\{ x \in [0,1]^{\mathbb{N}} : x_j = a \mbox{ for } j \in \mathbb{N} \setminus v \right\}
\end{displaymath}
for a finite set $\emptyset \neq v \subset \mathbb{N}$ and $a \in [0,1]$. We remind teh reader that essentially one only specifies those co-ordinates included in $v$, the remaining ones are specified via teh anchor point $a$.

{\em Variable subspace sampling} generalizes this idea to a sequence of finite-dimensional affine subspaces
\begin{displaymath}
\mathcal{X}_{v_1,a} \subset \mathcal{X}_{v_2,a} \subset \dots \, ,
\end{displaymath}
where $v=(v_i)_{i \in \mathbb{N}}$ is a given increasing sequence, $v_i \subset \mathbb{N}$, and $a \in [0,1]$. This sampling schemes allows us to choose integration nodes from subspaces of different dimensionality. To be able to compare introduce a cost function $c$, which si a mapping
\begin{displaymath}
c : [0,1]^{\mathbb{N}} \rightarrow \mathbb{N} \cap \left\{ \infty \right\}
\end{displaymath}
and which is to quantify how costly it is to evaluate teh integral $f$ at teh integration nodes. Following \cite{CDMGR09}, \cite{HMGNR10}, we formulate the cost of evaluating $f$ at the integration node $\bsx$ in terms of the dimension of teh finite-dimensional subspace from which the integration node is chosen. thsi means for fixed subspace sampling
\begin{equation} \label{eqcostevalfixed}
c_{v,a} = \left\{ \begin{array}{cc} \vert v \vert & , \mbox{ if } x \in \mathcal{X}_{v,a} \\ \infty & , \mbox{ otherwise.} \end{array} \right.
\end{equation}
For variable subspace sampling, which allows for the integration nodes to be chosen from a sequence of finite-dimensional subspaces, we choose the subspace with the smallest dimension in which the node lies,
\begin{equation} \label{eqcostevalvar}
c_{v,a}(x) = \inf \left\{ \dim (\mathcal{X}_{v_i,a} ) : x \in \mathcal{X}_{v_i,a} \right\} \ ,
\end{equation}
and set $\inf \emptyset = \infty$.

The randomized quadrature formulas employed in thsi note are based on scrambled polynomial lattice rules,
\begin{equation} \label{eqscrploylatt}
Q_{m,b,1:s}(f) = \frac{1}{b^m} \sum^{b^m}_{i=1} f(\bsx_i) \, ,
\end{equation}
where $\bsx_i \in [0,1)^s$ are obtained by scrambling a polynomial lattice rule, see \cite{BD10}. Defining
\begin{equation} \label{eqtruncf}
(\Psi_{v,a} f)(\bsx) = f(\bsx_v,a) \, ,
\end{equation}
we denote teh randomized quadrature formulas of interest in thsi note by
\begin{equation} \label{eqrandalg}
Q_{n,s,a} = Q_{\lfloor \log_b(n) \rfloor , b , 1:s} \circ \Psi_{1:s,a} \, ,
\end{equation}
where $n$ denotes the number of points of which the quadrature rule is comprised. Intuitively speaking, we carefully specify dimensions $1$ to $s$ via scrambled polynomial lattice rules and employ the anchor $a \in [0,1]$ for the subsequent dimensions. It is clear from the previous discussion, that for fixed subspace sampling, the notation introduced in \eqref{eqrandalg} is sufficient. For variable subspace sampling, however, we allow our integration nodes to be chosen from subspaces of different dimensions, in this sense the letter $s$ is not sufficient. We remark that this is addressed in Section \ref{secvarsubspacesampling}.

Next, we wish to discuss the cost of the randomized algorithms. As we only discuss randomized algorithms based on scrambled polynomial lattice rules in this note, we define the cost of fixed and variable subspace sampling under the premise that the randomized algorithm is based on a scrambled polynomial lattice rule. This simplifies the discussion, the cost model employed in \cite{HMGNR10}, which stems from \cite{CDMGR09}, allows for a much more general class of algorithms, see also \cite{KSWW09} for an even more general cost model.

Essentially, the cost of evaluating a randomized algorithm $Q$ is given by the sum of the costs of evaluating the function at the integration nodes chosen from the finite-dimensional subspaces. For teh fixed subspace sampling, assuming that nodes are chosen from $\mathcal{X}_{1:s,a}$
\begin{displaymath}
\mbox{cost}_{fix}(Q) = n s \, .
\end{displaymath}

For variable subspace sampling, we choose our integration nodes from a sequence of finite-dimensional affine subspaces, index $\bsv=(v_i)^{m}_{i=1}$, where $m \leq n$, as we use an $n$-point quadratuer rule, where $v_i=1:s_i$, $i=1,\dots,m$. For the subspace indexed vy $v_i$, eth integration nodes would be based on scrambled polynomial lattic erules whose integration nodes lie in $[0,1)^{s_i}$, and we denote the number of those integration nodes by $n_{v_i}$, where of course $\sum^m_{i=1} n_{v_i}=n$. Consequently, we have
\begin{equation} \label{eqcostalgvar}
\mbox{cost}_{var}(Q) = \sum^m_{i=1} s_i n_{v_i} \, .
\end{equation}

Finally, following \cite{HMGNR10}, for class of integrands $F$, where $f \in F$ is a mapping $f:\mathcal{X} \rightarrow \mathbb{R}$, $\mathcal{X} \subset [0,1]^{\mathbb{N}}$, we use the notation
\begin{displaymath}
I(f) = \int_{\mathcal{X}} f(\bsx) d\bsx
\end{displaymath}
and denote the worst-case error of a randomized algorithm $Q$, used to approximate integrands $f$ in the class $F$ by
\begin{displaymath}
e(Q,F) = \sup_{f \in F} \left( E \left( I(f) - Q(f) \right)^{2} \right)^{1/2} \, .
\end{displaymath}

Lastly, minimal errors, which are of great importance in information-based complexity, \cite{TWW88}, \cite{N88}, \cite{R00}, are defined by
\begin{displaymath}
e_{N, fix} (F) = \inf \left\{ e(Q,F) : \mbox{ cost }_{fix} (Q,F) \leq N  \right\}
\end{displaymath}
and
\begin{displaymath}
e_{N, var} (F) = \inf \left\{ e(Q,F) : \mbox{ cost }_{var}(Q,F) \leq N \right\} \, .
\end{displaymath}

The following result on numericla integration in $H(K_{1:s})$, see Equation \eqref{eqfindimHspace}, stems from \cite{BD10}.

\begin{theorem} \label{theoint} Assume $\sum^{\infty}_{j=1} \gamma^{\frac{1}{3 - \varepsilon}}_j < \infty$, for $\varepsilon >0$, then
\begin{displaymath}
e \left( Q_{b,m,1:s}, H(K_{1:s})  \right) \leq c_{\varepsilon} n^{-(3/2 - \varepsilon)} \, ,
\end{displaymath}
where $n=b^m$ and $Q_{b,m,1:s}$ is a scrambled polynomial lattice rule as defined in Equation \eqref{eqscrploylatt}.
\end{theorem}
\proof
From the proof of \cite[lemma 7]{YH05}, it is clear that the function space $H(K_{1:s})$ can be embedded in the space $V_{1,s,\gamma}$, as defined in \cite{BD10}, from which the result follows immediately.
\endproof

Concluding this section, we point out, that the forthcoming results will be presented for functions in $B(K))$, the unit ball in $H(K)$.

\section{Results on Fixed Subspace Sampling} \label{secfixsubspacesampl}

To fully specify the fixed subspace sampling algorithm, we only need to specify the dimension of the finite-dimensional subspace employed for sampling, and the number of integration nodes, which are based on a scrambled polynomial lattice rule, that we employ. As we wish to minimize worst-case errors for a fixed bound on the cost, say $N$, both, the dimension and the number of integration nodes are functions of $N$.

\begin{corollary}[Corollary 1, \cite{HMGNR10}]\label{corfixedsampling} Let $\varepsilon >0$ and let $\gamma_{j} \asymp j^{-\alpha}$, $\alpha \geq 3$. Choose
\begin{displaymath}
n \asymp N^{\frac{\alpha-1}{\alpha+2-\varepsilon}}
\end{displaymath}
and
\begin{displaymath}
s \asymp N^{\frac{3- \varepsilon}{\alpha+2-\varepsilon}}
\end{displaymath}
for $N \in \mathbb{N}$. Then, for $Q_N = Q_{n,s,a}$
\begin{displaymath}
e \left( Q_N, B(K) \right) \preceq N^{- \frac{(3 - \varepsilon)/2 (\alpha-1) }{\alpha +2 - \varepsilon}}
\end{displaymath}
and
\begin{displaymath}
\mbox{cost}_{fix} (Q_N, B(K)) \preceq N \, .
\end{displaymath}
\end{corollary}

\proof The result follows immediately from \cite[Theorem 1]{HMGNR10}, where we set $\alpha \geq 3$.
\endproof

\begin{remark}
In \cite{HMGNR10}, the same result was established under teh stronger assumption on the weights $\gamma_j \asymp j^{-\alpha}$, $\alpha > 4$. The result presented in Corollary \ref{corfixedsampling} is optimal for $\alpha \geq 3$, see \cite[Corollary 3]{HMGNR10}.
\end{remark}

\section{Results on Variable Subspace Sampling} \label{secvarsubspacesampling}

We carry out variable subspace sampling using the so-called multi-level approach, which was first introduced in \cite{H98}, \cite{H01}, see also \cite{G08a}, \cite{G08b}. teh idea underlying the multi-level approach is the following: We fix a sequence of sets
\begin{displaymath}
v_1 \subset \dots \subset v_L
\end{displaymath}
and the associated finite-dimensional affien subspaces
\begin{displaymath}
\mathcal{X}_{v_1} \subset \dots \subset \mathcal{X}_{v_L} \, .
\end{displaymath}
We use the integral associated with the finite-dimensional subspace of the largest dimension, $I( \Psi_{v_L,a} f)$ to approximate $I(f)$. However, we rewrite $I(\Psi_{v_L,a} f)$ as follows
\begin{displaymath}
I( \Psi_{v_L,a} f) = \sum^L_{l=1} I \left( \Psi_{v_l,a} f - \Psi_{v_{l-1},a} f \right)
\end{displaymath}
setting $\Psi_{v_0,a} f = 0$. Each of the integrals $I( \Psi_{v_l,a} f - \Psi_{v_{l-1},a} f )$ is now approximated using an independent randomized algorithm, based on a scrambled polynomial lattice rule, in particular we use a randomized algorithm
\begin{equation} \label{eqmultlevalg}
Q(f) = \sum^L_{l=1} Q_{n_l,s_l,a} ( f - \Psi_{1:s_{l-1} ,a} f ) \, ,
\end{equation}
so at level $l$, we use an algorithm based on a scrambled polynomial lattice rule consisting of $b^{\lfloor \log_b (n_l) \rfloor}$ points, which lie in $[0,1)^{s_l}$. The error associated with this algorithm can be split into bias and variance,
\begin{displaymath}
E \left( I(f) - Q(f) \right)^2 = \left( I(f) - I( \Psi_{1:s_L,a} f) \right)^2 + \Var(Q(f)) \, ,
\end{displaymath}
in particular
\begin{displaymath}
\Var(Q(f)) = \sum^L_{l=1} \Var \left( Q_{n_l,s_l,a} ( f - \Psi_{1:s_{l-1},a} f  ) \right) \, ,
\end{displaymath}
see \cite{HMGNR10}. Regarding the cost, from Equation \eqref{eqcostalgvar},
\begin{displaymath}
\mbox{ cost }_{var}(Q, B(K)) = \sum^L_{l=1} s_l n_l \, .
\end{displaymath}
By definition of variable subspace sampling, the dimension $s_l$ increases with $l$, but one would expect the variances $\Var \left( Q_{n_l,s_l,a} ( f - \Psi_{1:s_{l-1},a} f  ) \right)$ to decrease as $l$ increases; the challenge is to trade of these effects well.

\begin{corollary}[Corollary 4, \cite{HMGNR10}] \label{corvarsampling} Assume that $\gamma_j \asymp j^{- \alpha}$, for $\alpha > 3$, let $ 0 < \varepsilon < \min(6, \alpha -3)$ and put
\begin{displaymath}
\rho_1 = \frac{\alpha-1}{3 - \varepsilon/2} \, , \quad \rho_2 = \frac{\alpha -4 - \varepsilon}{3 - \varepsilon/2} \, .
\end{displaymath}
Choose $L$, $s_l$, $n_l$ according to Equations (26), (27), and (28) in \cite{HMGNR10}, and let $a \in [0,1]$. Take the corresponding multi-level algorithm $Q_N$ according to Equation \eqref{eqmultlevalg} based on the scrambled polynomial lattice rule. Then
\begin{displaymath}
e(Q_N, B(K)) \preceq \left\{ \begin{array}{cc} N^{-(3 - \varepsilon)/2} , & \mbox{ if } \alpha \geq 10 \, , \\ N^{-(3 - \varepsilon)/2 \frac{\alpha -1}{9}} \, , & \mbox{ if } \alpha < 10 \, , \end{array} \right.
\end{displaymath}
and
\begin{displaymath}
\mbox{cost}_{var}(Q_N, B(K)) \preceq N \, .
\end{displaymath}
\end{corollary}

\proof The proof follows immediately from \cite[Theorem 4]{HMGNR10}, with $\alpha' = 3 + \varepsilon$.
\endproof

\begin{remark} The same error bounds were established in \cite{HMGNR10}, but the rate $N^{(3 - \varepsilon)/2}$ was only established for $\alpha \geq 11$, whereas here it is achieved for $\alpha \geq 10$, due to an improved strong tractability result. Of course, for $\alpha \geq 10$, this result is optimal. Furthermore, we conclude that for (at least) $\alpha >7$, variable subspace sampling improves on fixed subspace sampling.
\end{remark}

\begin{remark} \label{remworstcasesett} We alert the reader to \cite{NHMGR10}, where infinite-dimensional integration in the worst-case setting is studied. In \cite{NHMGR10}, rank-$1$ lattice rules are employed as a basis for the algorithms, and we remark that in the worst-case setting, polynomial lattice rules have not been shown to improve on rank-1 lattice rules.
\end{remark}

\section{Future Research and Applications}

Applications of infinite-dimensional quadrature problems abound, to name but two from mathematical finance that received attention at MCQMC 2010, we list series expansions, such as the Karhunen-Lo\`eve expansion in the case of Brownian motion, see e.g. \cite{NH08} for a related publication, and the shot-noise representation in the context of L\'evy processes, see \cite{IK10}.

In the case of infinite-dimensional quadrature, the theory, as presented in \cite{CDMGR09}, \cite{KSWW09}, \cite{HMGNR10}, \cite{NHMGR10}, is ahead of the applications in the following sense: If we can specify the decay of the weights, all parameters of our algorithm are readily determined. However, in practice, the study of the weights has still not received enough attention, but it is clear that in order to fully reap the benefits of the algorithm, this needs to be addressed.

\end{document}